\documentclass{amsart}
\usepackage{amssymb,amsmath,latexsym,epsfig}

\newtheorem{thm}{Theorem}[section]
\newtheorem{dfn}[thm]{Definition}
\newtheorem{cor}[thm]{Corollary}
\newtheorem{prop}[thm]{Proposition}
\newtheorem{lemma}[thm]{Lemma}

\newcommand{\del}{\backslash}
\newcommand{\cl}{\hbox{\rm cl}}
\newcommand{\join}{\lor}
\newcommand{\meet}{\land}
\newcommand{\mc}{\mathcal}

\newfont{\menutt}{cmtt8}

\title[Transversal Lattices]{Transversal Lattices}
                   
\date{\today}
\author[J.~Bonin]
       {Joseph E.~Bonin}
\address
{Department of Mathematics\\ The George Washington University\\
Washington, D.C. 20052, USA} \email{jbonin@gwu.edu}

\begin{document}

\begin{abstract}
  A flat of a matroid is cyclic if it is a union of circuits; such
  flats form a lattice under inclusion and, up to isomorphism, all
  lattices can be obtained this way.  A lattice is a Tr-lattice if all
  matroids whose lattices of cyclic flats are isomorphic to it are
  transversal.  We investigate some sufficient conditions for a
  lattice to be a Tr-lattice; a corollary is that distributive
  lattices of dimension at most two are Tr-lattices.  We give a
  necessary condition: each element in a Tr-lattice has at most two
  covers.  We also give constructions that produce new Tr-lattices
  from known Tr-lattices.
\end{abstract}

\maketitle

\section{Introduction}

A flat $X$ of a matroid $M$ is \emph{cyclic} if the restriction $M|X$
has no isthmuses.  Ordered by inclusion, the cyclic flats form a
lattice, which we denote by $\mc{Z}(M)$.  Every lattice is isomorphic
to the lattice of cyclic flats of some (bi-transversal)
matroid~\cite{cyclic,julie}.  (All lattices considered in this paper
are finite.)  For certain lattices $L$, it is shown
in~\cite{acketa,ac} that if $\mc{Z}(M)$ is isomorphic to $L$, then the
matroid $M$ is transversal; lattices with this property are
\emph{transversal lattices} or \emph{Tr-lattices}.  In~\cite{cyclic},
lattices of width at most two are shown to be Tr-lattices.  In this
paper we treat a more general sufficient condition for a lattice to be
a Tr-lattice, we prove a necessary condition, and we show that the
class of Tr-lattices is closed under certain lattice operations.

Following a section of background, Section~\ref{sec:sufficient}
introduces MI-lattices and shows they are Tr-lattices.  Special cases
(e.g., distributive lattices of dimension at most two) are also
treated.  Section~\ref{sec:necessary} shows that each element of a
Tr-lattice has at most two covers.  Section~\ref{sec:examples} gives
ways to construct new MI-lattices (resp., Tr-lattices) from known
MI-lattices (resp., Tr-lattices).  Some open problems suggested by
this work are mentioned in the concluding section.

\section{Background}

We assume familiarity with basic matroid theory.  Our notation and
terminology for matroid theory follow~\cite{ox}; for ordered sets we
mostly follow~\cite{tomt}.  For a collection $\mc{F}$ of sets, we
write $\bigcap(\mc{F})$ for the intersection $\bigcap_{X\in \mc{F}}X$
and $\bigcup(\mc{F})$ for $\bigcup_{X\in \mc{F}}X$.  

Recall that every ordered set $P$ can be embedded in a product of
chains; the \emph{dimension} of $P$ is the least number of chains for
which there is such an embedding.  The lattices of dimension $2$ are
the \emph{planar} lattices: their Hasse diagrams can be drawn in the
plane without crossings (see, e.g.,~\cite[Chapter~3,
Theorem~5.1]{tomt}).  An \emph{antichain} in an ordered set is a
collection of mutually incomparable elements.  The \emph{width} of an
ordered set is the maximal cardinality among its antichains.  We say
$y$ is a \emph{cover} of $x$ in an ordered set $P$ if $x<y$ and there
is no $z$ in $P$ with $x<z<y$.  The least and greatest elements in a
lattice are denoted $\hat{0}$ and $\hat{1}$, respectively.  The
\emph{atoms} of a lattice are the elements that cover $\hat{0}$;
dually, the \emph{coatoms} are the elements that $\hat{1}$ covers. An
\emph{ideal} in an ordered set $P$ is a subset $I$ of $P$ such that if
$x\in I$ and $y\leq x$, then $y\in I$.  Dually, a \emph{filter} in $P$
is a subset $F$ such that if $x\in F$ and $y\geq x$, then $y\in F$.

It is well known and easy to see that while nonisomorphic matroids can
have the same cyclic flats, a matroid on a given set is determined by
its collection of cyclic flats along with their ranks.  In some cases
we will want to ignore the cyclic flats and instead focus on the ranks
assigned to the elements of an abstract lattice; this is justified by
the following special case of~\cite[Theorem 1]{julie}.

\begin{prop}\label{prop:sims}
  Let $L$ be a lattice.  Given $\rho:L\rightarrow \mathbb{Z}$ with
  \begin{itemize}
  \item[(a)] $\rho(\hat{0})=0$,
  \item[(b)] $\rho(x)<\rho(y)$ whenever $x<y$, and
  \item[(c)] $\rho(x\join y) + \rho(x\meet y)\leq \rho(x) + \rho(y)$
    whenever $x$ and $y$ are incomparable,
  \end{itemize}
  there is a matroid $M$ and an isomorphism
  $\phi:L\rightarrow\mc{Z}(M)$ with $\rho(x)=r\bigl(\phi(x)\bigr)$.
\end{prop}

A key result we use to prove that certain lattices are (or are not)
Tr-lattices is the following characterization of transversal matroids,
which was first formulated by Mason using cyclic sets and later
refined to cyclic flats by Ingleton~\cite{ingleton}.  (The statement
in~\cite{ingleton} uses all nonempty collections of cyclic flats, but
an elementary argument shows that it suffices to consider nonempty
antichains of cyclic flats; see the discussion
after~\cite[Lemma~5.6]{cyclic}.)

\begin{prop}\label{prop:mi} 
  A matroid $M$ is transversal if and only if for every nonempty
  antichain $\mc{A}$ in $\mc{Z}(M)$,
    \begin{equation}\label{mi}
      r\bigl(\bigcap(\mc{A})\bigr) \leq \sum_{\mc{F}\subseteq\mc{A}}
      (-1)^{|\mc{F}|+1} r\bigl(\bigcup(\mc{F})\bigr).\tag{MI}
    \end{equation}
\end{prop}

The join in $\mc{Z}(M)$ (as in the lattice of flats) is given by
$A\join B=\cl(A\cup B)$, so one can replace the alternating sum in
inequality~(MI) by the corresponding alternating sum of ranks of joins
of cyclic flats. 

Unlike in the lattice of flats, the meet operation in $\mc{Z}(M)$
might not be intersection: $X\meet Y$ is the union of the circuits
that are contained in $X\cap Y$.

Since the complements of the flats of a matroid are the unions of its
cocircuits, $X$ is a cyclic flat of $M$ if and only if $E(M)-X$ is a
cyclic flat of the dual, $M^*$.  Thus, $\mc{Z}(M^*)$ is isomorphic to
the order dual of $\mc{Z}(M)$.

Let $S$ and $E$ be the least and greatest cyclic flats of $M$.  Note
that for $X\in \mc{Z}(M)$, the lattice $\mc{Z}(M|X)$ is the interval
$[S,X]$ in $\mc{Z}(M)$ and, dually, the lattice $\mc{Z}(M/X)$ is
isomorphic to the interval $[X,E]$ in $\mc{Z}(M)$ via the isomorphism
$Y\mapsto Y\cup X$.  (The lattices of cyclic flats of other minors are
not as simple to describe.)

\section{Sufficient conditions for a lattice to be a
  Tr-lattices}\label{sec:sufficient}

To convey the spirit of the main result of this section
(Theorem~\ref{thm:primitive}) before defining the technical condition
involved, we cite the following theorem, which, as we will show, is
implied by the main result.

\begin{thm}\label{thm:sublattice}
  If (a) $\mc{Z}(M)$ has dimension at most two and (b) for each
  antichain $\mc{A}$ of $\mc{Z}(M)$, the sublattice of $\mc{Z}(M)$
  generated by $\mc{A}$ is distributive, then $M$ and all of its
  minors, as well as their duals, are transversal.
\end{thm}

\begin{cor}
  If $\mc{Z}(M)$ is distributive and has dimension at most two, then
  $M$ and all of its minors, as well as their duals, are transversal.
\end{cor}

The main result of this section uses the following notions.

\begin{dfn}\label{def:mi}
  An \emph{MI-ordering} of an antichain $\mc{A}$ in a lattice $L$ is a
  permutation $a_1,a_2,\ldots,a_t$ of $\mc{A}$ so that
  \begin{itemize}
  \item[(i)] $a_i\join a_{i+1}\join \cdots\join a_k = a_i\join a_k$
    for $1\leq i<k\leq t$ and
  \item[(ii)] $(a_1\meet a_2\meet\cdots\meet a_k)\join
    a_{k+1}=a_k\join a_{k+1}$ for $1<k<t$.
  \end{itemize}
  An antichain is \emph{MI-orderable} if it has an MI-ordering.  A
  lattice is \emph{MI-orderable}, or is an \emph{MI-lattice}, if each
  of its antichains is MI-orderable.
\end{dfn}

\begin{thm}\label{thm:primitive}
  Let $M$ be a matroid.
  \begin{enumerate}
  \item[(i)] Each MI-orderable antichain in $\mc{Z}(M)$ satisfies
    inequality~(MI).
  \item[(ii)] If $\mc{Z}(M)$ is MI-orderable, then $M$ and all of its
    minors are transversal.
  \end{enumerate}
\end{thm}

\begin{cor}
  MI-lattices are Tr-lattices.
\end{cor}

\begin{figure}[t]
\begin{center}
\epsfxsize 10cm \epsffile{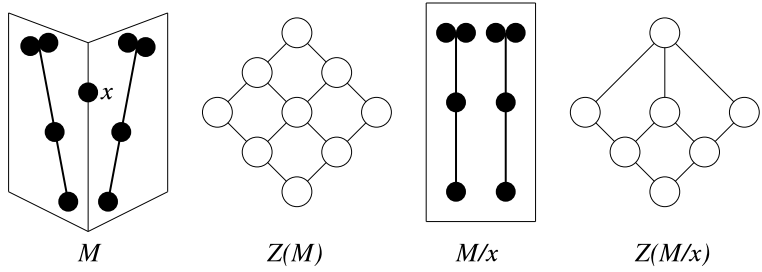} 
\caption{The lattice of cyclic flats of a matroid $M$ and that of
  $M/x$.}\label{fig:minor}
\end{center}
\end{figure}

Before proving Theorem~\ref{thm:primitive}, we note a subtlety that
explains the approach we take to prove part (ii): if $N$ is a minor of
$M$ and $\mc{Z}(M)$ is MI-orderable, then $\mc{Z}(N)$ may or may not
be MI-orderable.  Indeed, $\mc{Z}(N)$ may not even be a Tr-lattice,
and this is so even for deletions of $M$.  (Recall that the class of
transversal matroids is closed under deletions but not under
contractions, so one might expect deletions to be somewhat more tame.)
For example, for the matroid $M$ in Figure~\ref{fig:minor},
$\mc{Z}(M)$ is MI-orderable.  Since this lattice is isomorphic to its
order dual, $\mc{Z}(M^*)$ is also MI-orderable.  The lattice
$\mc{Z}(M/x)$ is also shown; by checking directly or applying
Theorem~\ref{thm:ideals}, we have that $\mc{Z}(M/x)$ is MI-orderable.
However, by Theorem~\ref{thm:cover}, its order dual, which is
$\mc{Z}(M^*\del x)$, is not a Tr-lattice.  This example also shows
that the minor-closed, dual-closed class of matroids described in
Theorem~\ref{thm:sublattice} is not determined by lattice-theoretic
properties that apply to the lattices of cyclic flats of all matroids
in the class.

We prove Theorem~\ref{thm:primitive} via a sequence of lemmas.  The
first lemma gives a rank inequality associated with each MI-orderable
antichain of $\mc{Z}(M)$.  Note that for two-element antichains, this
inequality is the semimodular inequality.  (The meet and join
operations in this and other results are in $\mc{Z}(M)$.)

\begin{lemma}\label{lemma:in}
  Let $A_1,A_2,\ldots,A_t$ be an antichain of cyclic flats in a
  matroid $M$ such that $(A_1\meet A_2\meet\cdots\meet A_k)\join
  A_{k+1}=A_k\join A_{k+1}$ whenever $1\leq k<t$.  Then for $k$ with $k\leq t$,
  \begin{equation}\label{dkeyineq} r(A_1\cap A_2\cap \cdots
    \cap A_k)\leq \sum_{i=1}^k r(A_i)-\sum_{i=1}^{k-1}r(A_i\cup
    A_{i+1}).
  \end{equation}
\end{lemma}

\begin{proof}
  We prove the inequality by induction on $k$.  Equality holds for
  $k=1$. Assume the result holds in case $k$.  Semimodularity gives
  $$r(A_1\cap A_2\cap \cdots \cap A_{k+1})+ r\bigl((A_1\cap A_2\cap
  \cdots \cap A_k)\cup A_{k+1}\bigr)- r(A_1\cap A_2\cap \cdots \cap
  A_k)\leq r(A_{k+1}).$$
  Adding this inequality to
  inequality~(\ref{dkeyineq}) gives $$r(A_1\cap A_2\cap \cdots \cap
  A_{k+1}) + r\bigl((A_1\cap A_2\cap \cdots \cap A_k)\cup
  A_{k+1}\bigr)\leq \sum_{i=1}^{k+1} r(A_i)-\sum_{i=1}^{k-1}r(A_i\cup
  A_{i+1}),$$
  so if we show $r\bigl((A_1\cap A_2\cap \cdots \cap
  A_k)\cup A_{k+1}\bigr)=r(A_k\cup A_{k+1})$, then the inequality we
  want follows.  This equality holds since $A_1\meet A_2\meet \cdots
  \meet A_k\subseteq A_1\cap A_2\cap \cdots \cap A_k \subseteq A_k$
  and $(A_1\meet A_2\meet \cdots \meet A_k)\join A_{k+1}=A_k\join
  A_{k+1}$.
\end{proof}

\begin{lemma}\label{lemma:t}
  If an antichain $\mc{A}$ in $\mc{Z}(M)$ can be ordered as
  $A_1,A_2,\ldots,A_t$ so that
  \begin{itemize}
  \item[(i)] $A_i\join A_{i+1}\join \cdots\join A_k = A_i\join A_k$
    whenever $1\leq i<k\leq t$ and
  \item[(ii)] $r(A_1\cap A_2\cap \cdots \cap A_t) \leq \sum_{i=1}^t
    r(A_i)-\sum_{i=1}^{t-1}r(A_i\cup A_{i+1}),$
  \end{itemize}
  then $\mc{A}$ satisfies inequality~(MI).
\end{lemma}

\begin{proof}
  Assume properties (i) and (ii) hold.  For $1\leq i\leq j\leq t$, set
  $$\mc{A}_{i,j}= \{\mc{F}\,:\,\mc{F}\subseteq\mc{A},\,\,
  i=\min(k\,:\,A_k\in\mc{F}), \text{ and } j=\max(k\,:\,A_k\in\mc{F})
  \}.$$
  Thus, if $\mc{F}\in \mc{A}_{i,j}$, then
  $\cl\bigl(\bigcup(\mc{F})\bigr)=A_i\join A_j$.  If $j>i+1$, then the
  terms on the right side of inequality~(MI) that arise from the sets
  in $\mc{A}_{i,j}$ cancel since there is a parity-switching
  involution $\phi$ of $\mc{A}_{i,j}$: fix $k$ with $i<k<j$ and let
  $$\phi(\mc{F})= \left\{
    \begin{array}{ll} \mc{F}\cup\{A_k\},
      &\mbox{if $A_k\not\in\mc{F}$;}\\
      \mc{F}-\{A_k\}, &\mbox{if $A_k\in\mc{F}$.}
    \end{array}\right.$$
  Thus, inequality~(MI) reduces to the inequality that is assumed in
  property (ii).
\end{proof}

The previous two lemmas show that MI-lattices are Tr-lattices.  To
prove the stronger assertion in part (ii) of
Theorem~\ref{thm:primitive}, we show that if the antichains in
$\mc{Z}(M)$ satisfy the hypotheses of Lemma~\ref{lemma:t}, then so do
the antichains of single-element deletions and single-element
contractions of $M$.  (Note that unlike the hypotheses of
Theorem~\ref{thm:primitive}, condition (ii) in Lemma~\ref{lemma:t} is
not a lattice-theoretic property.)  We start with a lemma about the
cyclic flats of such minors.

\begin{lemma}~\label{lem:cydel}
  For an element $x$ of $M$ and a cyclic flat $A$ of either $M\del x$
  or $M/x$, the flat $\bar{A}=\cl_M(A)$ of $M$ is cyclic;
  furthermore, $\bar{A}$ is either $A$ or $A\cup x$, so $\bar{A}-x=A$.
\end{lemma}

\begin{proof}
  For a cyclic flat $A$ of $M\del x$, the assertions are transparent.
  Let $A$ be a cyclic flat of $M/x$ and let $S$ be the ground set of
  $M/x$.  Thus, $S-A$ is a cyclic flat of the dual of $M/x$, that is,
  of $M^*\del x$, so $\cl_{M^*}(S-A)$, which is either $S-A$ or
  $(S-A)\cup x$, is a cyclic flat of $M^*$.  Therefore either $A\cup
  x$ or $A$ is a cyclic flat of $M$, from which the result follows.
\end{proof}

\begin{lemma}
  If each antichain in $\mc{Z}(M)$ can be ordered so that properties
  \emph{(i)} and \emph{(ii)} of Lemma~\emph{\ref{lemma:t}} hold, then
  the same is true for each antichain in $\mc{Z}(M\del x)$ and each
  antichain in $\mc{Z}(M/x)$.
\end{lemma}

\begin{proof}
  The proofs for $\mc{Z}(M\del x)$ and $\mc{Z}(M/x)$ are similar and,
  since each deletion of a transversal matroid is transversal, only
  the result about contractions is needed to prove
  Theorem~\ref{thm:primitive}, so we treat only $\mc{Z}(M/x)$.  We use
  the notation $\bar{A}$ of Lemma~\ref{lem:cydel}.  Let $\mc{A}$ be an
  antichain in $\mc{Z}(M/x)$.  Note that $\{\bar{A}\,:\,A\in\mc{A}\}$
  is an antichain in $\mc{Z}(M)$.  By hypothesis, there is an ordering
  $A_1,A_2,\ldots,A_t$ of $\mc{A}$ so that in $M$ and $\mc{Z}(M)$,
  \begin{equation}\label{del1}
    \bar{A}_i\join \bar{A}_{i+1}\join \cdots\join \bar{A}_k =
    \bar{A}_i\join \bar{A}_k, \quad\text{ for }  1\leq i<k\leq t,
  \end{equation}
  and 
  \begin{equation}\label{del2}
    r_M(\bar{A}_1\cap \bar{A}_2\cap
    \cdots \cap \bar{A}_t)+\sum_{i=1}^{t-1}r_M(\bar{A}_i\cup
    \bar{A}_{i+1})\leq \sum_{i=1}^t r_M(\bar{A}_i). 
  \end{equation}
  Since $\bar{A}_j=\cl_M(A_j)$ and since $A\join B$ in $\mc{Z}(M)$ is
  $\cl_M(A\cup B)$, by equation~(\ref{del1}) $A_i\cup
  A_{i+1}\cup\cdots\cup A_k$ and $A_i\cup A_k$ have the same closure
  in $M$, and so in $M/x$; thus, as needed, $A_i\join A_{i+1}\join
  \cdots\join A_k = A_i\join A_k$ in $\mc{Z}(M/x)$.  The rank
  inequality in $M/x$ is immediate if $x$ is a loop of $M$, so assume
  this is not the case.  Assume $x$ is in exactly $h$ of the cyclic
  flats $\bar{A}_1,\bar{A}_2,\ldots,\bar{A}_t$ of $M$.  Thus,
  $$h+\sum_{i=1}^t r_{M/x}(A_i)=\sum_{i=1}^t r_M(\bar{A}_i).$$
  Also,
  since $x$ is in at least $h$ of the sets $\bar{A}_1\cap
  \bar{A}_2\cap \cdots \cap \bar{A}_t$ and $\bar{A}_i\cup
  \bar{A}_{i+1}$, we have
  \begin{align}
    h+ r_{M/x}(A_1\cap A_2\cap& \cdots \cap A_t)
    +\sum_{i=1}^{t-1}r_{M/x}(A_i\cup A_{i+1}) \notag  \\ & \leq
    r_M(\bar{A}_1\cap \bar{A}_2\cap \cdots \cap
    \bar{A}_t)+\sum_{i=1}^{t-1}r_M(\bar{A}_i\cup \bar{A}_{i+1}).
    \notag
  \end{align}
  The last two conclusions and inequality (\ref{del2}) give
  $$r_{M/x}(A_1\cap A_2\cap \cdots \cap A_t)
  +\sum_{i=1}^{t-1}r_{M/x}(A_i\cup A_{i+1})\leq \sum_{i=1}^t
  r_{M/x}(A_i),$$
  as needed.
\end{proof}

The lemmas above complete the proof of Theorem~\ref{thm:primitive}.
We now show that Theorem~\ref{thm:sublattice} follows.  Recall that
$\mc{Z}(M^*)$ is isomorphic to the order dual of $\mc{Z}(M)$, so $M^*$
satisfies the hypotheses of Theorem~\ref{thm:sublattice} if and only
if $M$ does.  Thus, the next lemma suffices to prove
Theorem~\ref{thm:sublattice}.

\begin{lemma}\label{lemma:o}
  If a lattice $L$ has dimension at most two and each of its
  antichains generates a distributive sublattice, then $L$ is MI-orderable.
\end{lemma}

\begin{proof}
  We may assume $L$ is a suborder of $\mathbb{N}^2$.  List the
  elements of an antichain as $a_1,a_2,\ldots,a_t$ where
  $a_i=(x_i,y_i)$ with $x_1>x_2>\cdots>x_t$; thus,
  $y_1<y_2<\cdots<y_t$.  Clearly $a_i\join a_{i+1}\join \cdots\join
  a_k \geq a_i\join a_k$.  Let $a_i\join a_k$ be $(p,q)$.  Thus,
  $p\geq x_i$ and $q\geq y_k$, so $(p,q)\geq (x_j,y_j)$ for $i\leq
  j\leq k$, and so $a_i\join a_{i+1}\join \cdots\join a_k = a_i\join
  a_k$.  One gets $a_i\meet a_{i+1}\meet \cdots\meet a_k = a_i\meet
  a_k$ similarly, so property (ii) of Definition~\ref{def:mi} can be
  rewritten as $(a_1\meet a_k)\join a_{k+1} = a_k\join a_{k+1}$, or,
  by the distributive law, $(a_1\join a_{k+1})\meet (a_k\join
  a_{k+1})=a_k\join a_{k+1}$, that is, $a_1\join a_{k+1} \geq a_k\join
  a_{k+1}$.  This property holds since $a_1\join a_{k+1}=a_1\join
  a_2\join \cdots\join a_k\join a_{k+1}\geq a_k\join a_{k+1}$.
\end{proof}

Antichains of at most two elements are trivially MI-orderable, so
Theorem~\ref{thm:primitive} has the following corollary (as do
Theorem~\ref{thm:sublattice} and Proposition~\ref{prop:mi}).

\begin{cor}\label{cor:width2} \emph{\cite[Theorem
    5.7.]{cyclic}}  Lattices of width at most two are Tr-lattices.
\end{cor}

Section~\ref{sec:examples} includes examples of lattices to which
Theorem~\ref{thm:primitive} but not Theorem~\ref{thm:sublattice}
applies, as well as Tr-lattices that are not MI-lattices.

\section{A necessary condition for a lattice to be a
  Tr-lattice}\label{sec:necessary}

Condition (ii) of Definition~\ref{def:mi} is violated by any three
covers of a given element, so each element of an MI-lattice has at
most two covers.  In this section, we show that the same is true of
any Tr-lattice.  (The examples in the next section show there is no
bound on the number of elements that an element in a Tr-lattice
covers.)

\begin{thm}\label{thm:cover}
  Each element of a Tr-lattice has at most two covers.
\end{thm}

\begin{proof}
  Let the element $x$ of a lattice $L$ have at least three covers.  We
  prove that $L$ is not a Tr-lattice by defining a function
  $\rho:L\rightarrow \mathbb{Z}$ so that properties (a)--(c) in
  Proposition~\ref{prop:sims} hold and inequality (MI) fails.  For
  $y\in L$, let $F_y$ be the principal filter $\{u\,:u\geq y\}$ in
  $L$.  Thus, the sublattice $F_x$ of $L$ has at least three atoms.

  Define $\rho':L\rightarrow \mathbb{Z}$ by $\rho'(y) =
  \bigl|L-F_y\bigr|$.  It follows easily that $\rho'$ satisfies
  properties (a)--(c) in Proposition~\ref{prop:sims}.  For $u,v,w\in
  F_x-\{x\}$, let $$m(u,v,w)=\rho'(u)+\rho'(v)+\rho'(w)-\rho'(u\join
  v)-\rho'(u\join w)-\rho'(v\join w)+\rho'(u\join v\join
  w)-\rho'(x).$$
  By inclusion-exclusion, $m(u,v,w)=|F_x-(F_u\cup
  F_v\cup F_w)|$.  Set $$k=\,\min\{m(u,v,w)\,:\,u,v,w > x\} =
  |F_x|-\max\{\bigl|F_u\cup F_v\cup F_w\bigr| \,:\,u,v,w > x\}.$$
  Thus, $k$ is the minimal size of the complement, in $F_x$, of the
  union of three proper principal filters in $F_x$.  Note that if
  $k=m(u,v,w)$, then $u,v,w$ are distinct covers of $x$.  Define
  $\rho:L\rightarrow\mathbb{Z}$ by
  \begin{eqnarray}
    \notag
    \rho(y) = 
    \left\{
      \begin{array}{ll}
        \rho'(y), &\mbox{if $y\leq x$,}\\
        \rho'(x)-k-1, &\mbox{otherwise}.
      \end{array}
    \right.
  \end{eqnarray}
  Clearly $\rho$ satisfies property (a) of
  Proposition~\ref{prop:sims}.  Properties (b) and (c) for $\rho$
  follow from these properties for $\rho'$ except in two cases, which
  we address below:
  \begin{itemize}
  \item[(i)] $\rho(y)<\rho(z)$ if $y<z$, $y\leq x$, and $z\not\leq x$,
    and
  \item[(ii)] $\rho(y)+\rho(z)\geq \rho(y\join z)+ \rho(y\meet z)$ if
    $y\not\leq x$, $z\not\leq x$, and $y\meet z\leq x$.
  \end{itemize} 
  
  Assume $y<z$, $y\leq x$, and $z\not\leq x$.  Thus, $F_x\subseteq
  F_y$.  The inequality in statement (i) reduces to $k+2\leq \rho'(z)
  - \rho'(y) = |F_y-F_z|$.  Note that $F_z\cap F_x$ is the principal
  filter $F_{z\join x}$, which, since $z\not \leq x$, is properly
  contained in $F_x$; thus, there are at least $k+2$ elements in
  $F_x-F_z$, and so in $F_y-F_z$, which proves statement (i).

  Now assume $y\not\leq x$, $z\not\leq x$, and $y\meet z\leq x$.  The
  inequality in statement (ii) is $$|L|-|F_y|-k-1 + |L|-|F_z|-k-1\geq
  |L|-|F_{y\join z}|-k-1 + |L|-|F_{y\meet z}|,$$
  that is, $|F_{y\meet
    z}-(F_y\cup F_z)|\geq k+1.$ Note that $F_x\subseteq F_{y\meet z}$
  and $$(F_y\cup F_z)\cap F_x = (F_y\cap F_x)\cup(F_z\cap
  F_x)=F_{y\join x}\cup F_{z\join x},$$
  which is the union of two
  principal filters that are properly contained in $F_x$; thus, there
  are at least $k+1$ elements in $F_x-(F_y\cup F_z)$ and so in
  $F_{y\meet z}-(F_y\cup F_z)$, which proves statement (ii).

  Let $M$ be a matroid arising from $L$ and $\rho$ as in
  Proposition~\ref{prop:sims}.  Fix $u,v,w$ with $k=m(u,v,w)$ and let
  $U$, $V$, and $W$ be the corresponding cyclic flats of $M$.  The
  definitions of $m$ and $\rho$ give $$r(U)+r(V)+r(W)-r(U\cup
  V)-r(U\cup W)-r(V\cup W)+r(U\cup V\cup W)=r(X)-1.$$
  Since $r(X)\leq
  r(U\cap V\cap W)$, it follows that the antichain $\{U,V,W\}$ of
  $\mc{Z}(M)$ does not satisfy inequality (MI).  Thus, $M$ is not
  transversal, so $L$ is not a Tr-lattice.
\end{proof}

A matroid $M$ is \emph{nested} if $\mc{Z}(M)$ is a chain.  These
matroids have arisen many times in a variety of contexts
(see~\cite[Section 4]{lpm2} for more information).  That
$\mc{Z}(M\oplus N)$ is isomorphic to the product
$\mc{Z}(M)\times\mc{Z}(N)$ gives the following corollary.

\begin{cor}
  If $\mc{Z}(M)$ is a Tr-lattice, then the matroid obtained from $M$
  by deleting all loops and isthmuses is either a direct sum of at
  most two nested matroids or it is connected.
\end{cor}

\section{Examples and constructions}\label{sec:examples}

This section gives examples of MI-lattices to which
Theorem~\ref{thm:sublattice} does not apply and Tr-lattices that are
not MI-lattices.  We also show how to construct new MI-lattices from
given MI-lattices, and likewise for Tr-lattices.

\begin{figure}[t]
\begin{center}
\epsfxsize 9.2cm \epsffile{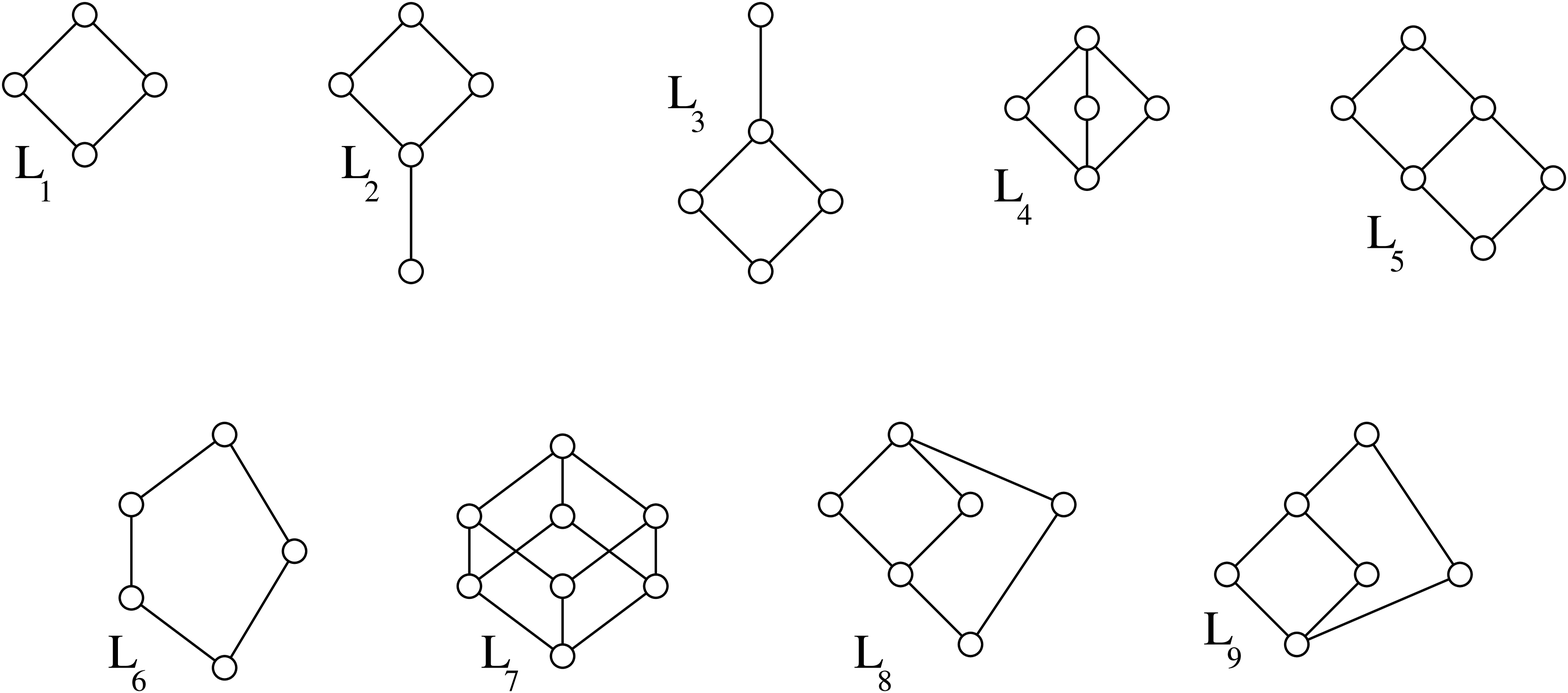} 
\caption{The lattices Acketa considered.}\label{fig:lattices}
\end{center}
\end{figure}

\begin{figure}[t]
\begin{center}
\epsfxsize 4.8cm \epsffile{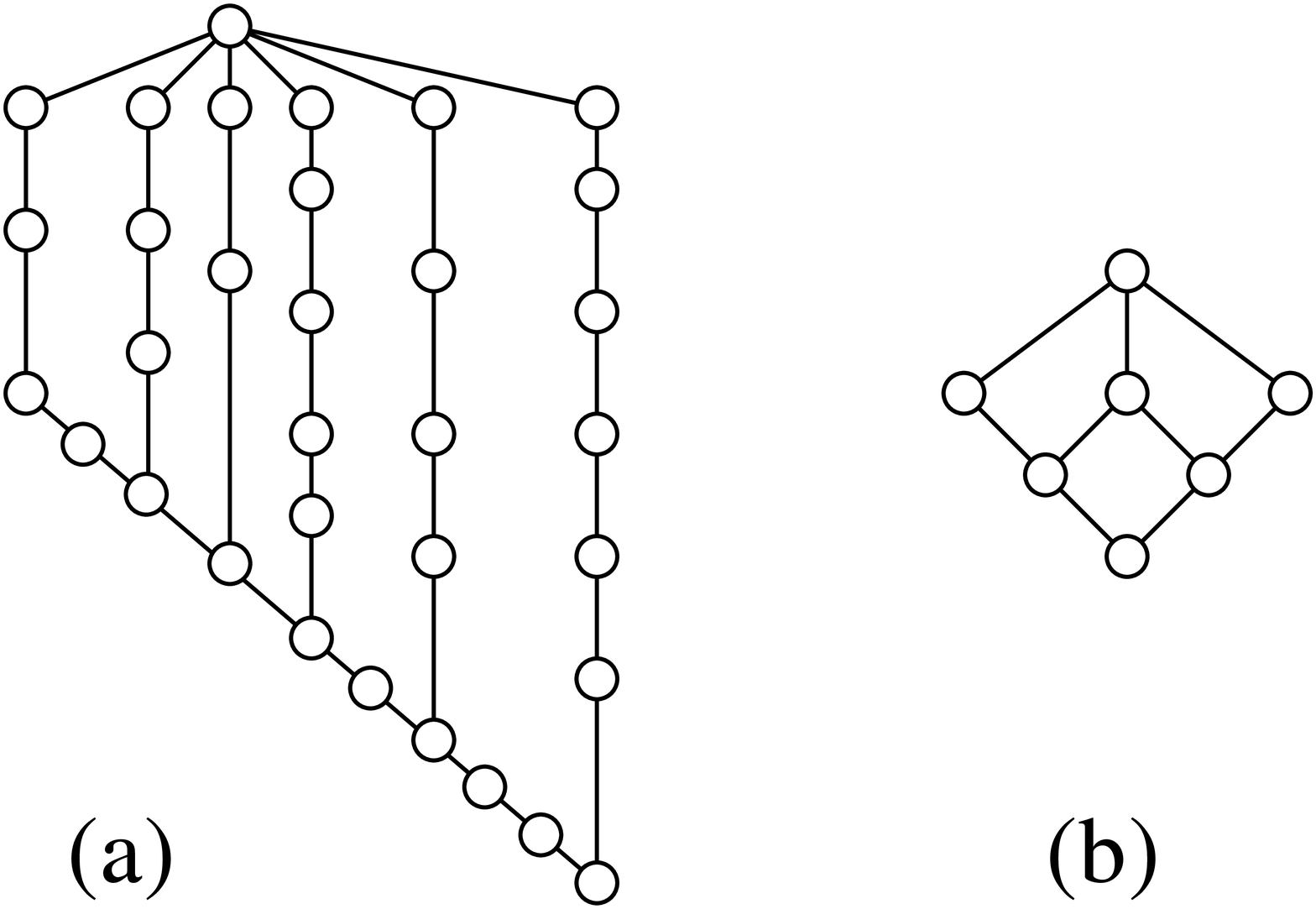} 
\caption{(a): A generic lattice like $L_8$. (b): A lattice $L_I$
  obtained from an ideal in a product of two three-element
  chains.}\label{fig:examples1}
\end{center}
\end{figure}

Acketa~\cite{acketa,ac} proved that chains and the lattices $L_1$,
$L_2$, and $L_3$ of Figure~\ref{fig:lattices} are Tr-lattices
(Corollary~\ref{cor:width2} applies); he noted that $L_4$ is not a
Tr-lattice; he conjectured that $L_5$, $L_6$, and $L_7$ are
Tr-lattices (Corollary~\ref{cor:width2} applies to $L_5$ and $L_6$;
Theorem~\ref{thm:cover} shows that $L_7$ is not a Tr-lattice); he
proved that $L_8$ is a Tr-lattice; he also showed that $L_9$ (the dual
of $L_8$) is not a Tr-lattice.  We note that $L_8$ is in an infinite
family of MI-lattices; Figure~\ref{fig:examples1}.a gives another such
lattice.  The defining properties of these lattices are that the
interval between $\hat{0}$ and any coatom is a chain, and for one of
these chains (e.g., the left-most chain in
Figure~\ref{fig:examples1}.a), all other such chains intersect it in
different initial segments.

Sublattices of MI-lattices are clearly MI-lattices.  The next result
gives another simple construction for MI-lattices.  (See
Figure~\ref{fig:examples1}.b.)

\begin{thm}\label{thm:ideals}
  For any ideal $I$ in an MI-lattice $L$, the lattice $L_I$ induced on
  the set $I\cup\{\hat{1}\}$ by the same order is MI-orderable.
\end{thm}

\begin{proof}
  Each antichain $\mc{A}$ of $L_I$ is an antichain of $L$; order
  $\mc{A}$ so that properties (i) and (ii) of Definition~\ref{def:mi}
  hold in $L$.  Let $z$ be the join of $\{a_i,a_{i+1},\ldots,a_k\}$
  and of $\{a_i,a_k\}$ in $L$.  If $z \in I$, then $z$ is the join of
  each of these sets in $L_I$, otherwise both sets have join $\hat{1}$
  in $L_I$.  Thus, property (i) holds in $L_I$.  The same ideas show
  that property (ii) holds in $L_I$ since the meet operations are
  identical in $L$ and $L_I$.
\end{proof}

Recall that the linear sum (or ordinal sum) of partial orders $P$ and
$Q$, where $P$ and $Q$ are disjoint, is the order on $P\cup Q$ in
which the restriction to $P$ is the order on $P$, the restriction to
$Q$ is the order on $Q$, and every element of $P$ is less than every
element of $Q$.  The following result is immediate.

\begin{thm}\label{thm:linsum}
  The class of MI-lattices is closed under linear sums.
\end{thm}

The same result holds for the closely-related operation that, given
lattices $L$ and $L'$, forms the Hasse diagram of the new lattice from
those of $L$ and $L'$ by identifying the greatest element of $L$ with
the least element of $L'$.  It follows from Theorem~\ref{thm:lex}
below that the same two results hold for Tr-lattices.  By
Theorem~\ref{thm:cover}, the class of MI-lattices and the class of
Tr-lattices are not closed under direct products.

\begin{figure}[t]
\begin{center}
\epsfxsize 8.5cm \epsffile{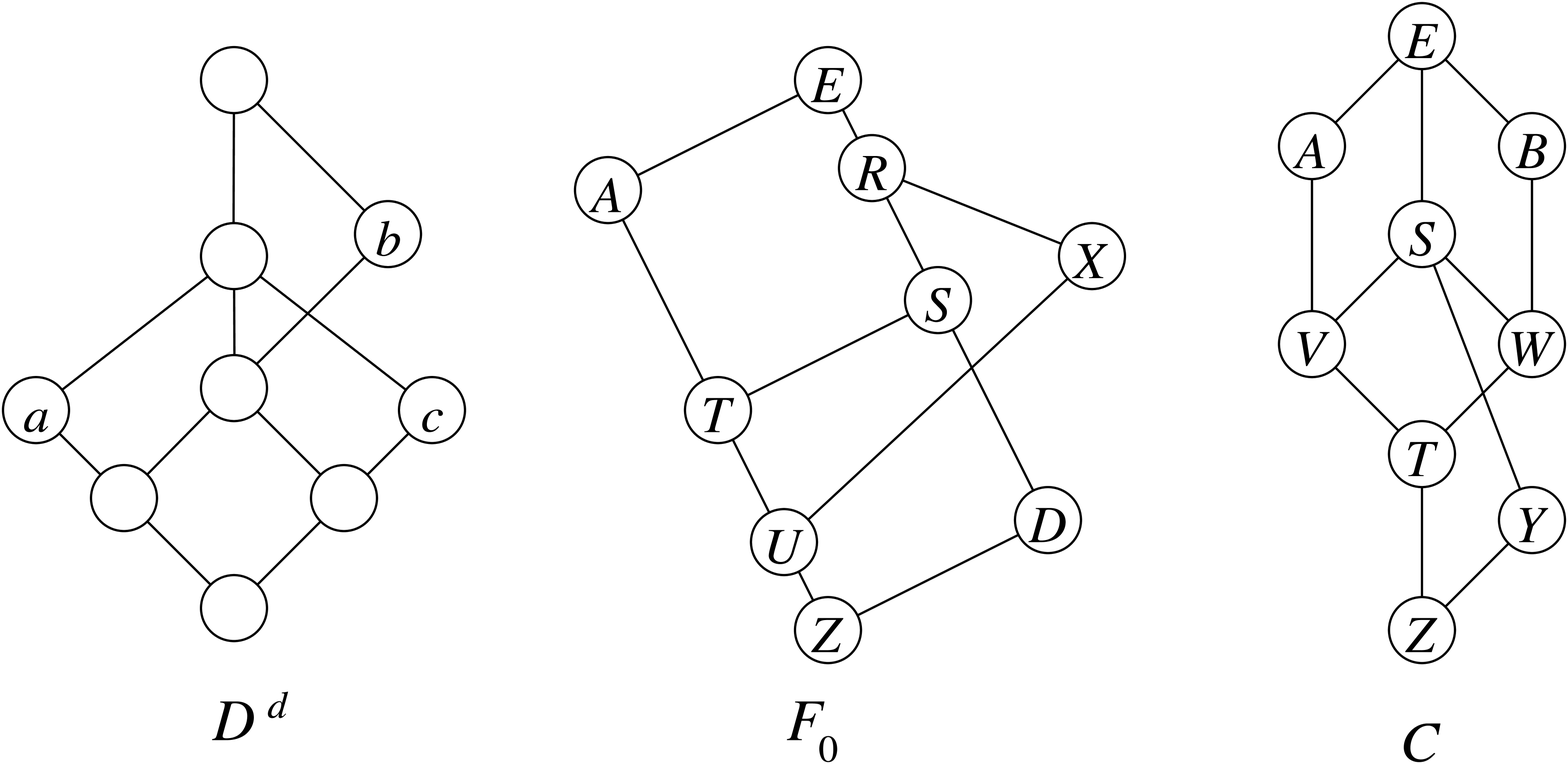} 
\caption{Three nonplanar Tr-lattices; only $D^d$ is
  MI-orderable.}\label{fig:examples2}
\end{center}
\end{figure}

We next treat three particular Tr-lattices of dimension $3$, only one
of which is MI-orderable.  These lattices, which are shown in
Figure~\ref{fig:examples2}, are among the forbidden sublattices for
planar lattices (see~\cite{kr}).  (No other forbidden sublattices for
planar lattices satisfy the necessary condition for Tr-lattices given
in Theorem~\ref{thm:cover}.)

\begin{thm}\label{thm:nonplanar}
  The lattice $D^d$ is MI-orderable.  The lattices $F_0$ and $C$ are
  Tr-lattices that are not MI-orderable.
\end{thm}

\begin{proof}
  The sublattice of $D^d$ formed by removing $b$ is the linear sum of
  the lattice in Figure~\ref{fig:examples1}.b and a single-element
  lattice.  Since this linear sum is MI-orderable, we need only check
  that the antichains in $D^d$ that contain $b$ are MI-orderable.
  Antichains of two elements are automatically MI-orderable; the only
  larger antichain in $D^d$ that includes $b$ is $\{a,b,c\}$, for
  which $b,a,c$ is an MI-ordering.

  In $F_0$, the antichains of more than two elements are $\{A,S,X\}$,
  $\{X,T,D\}$, and $\{X,A,D\}$.  The first two are MI-orderable
  (ordered as written), so we need only show that in any matroid $M$
  for which $\mc{Z}(M)$ is isomorphic to $F_0$, the flats
  corresponding to $X,A,D$ (for which we use the same notation)
  satisfy inequality (MI), which in this case is
  $r(X)+r(A)+r(D)-r(R)-r(E)\geq r(X\cap A\cap D)$.  By semimodularity,
  $$r(A)+r(S)\geq r(E)+r(A\cap S).$$
  The inclusions $T\subseteq A\cap
  S\subseteq S$ give $\cl\bigl((A\cap S)\cup X\bigr)=R$, so $$r(A\cap
  S)+r(X)\geq r(R)+r(A\cap S\cap X).$$
  The inclusions $U\subseteq
  A\cap S\cap X\subseteq S$ give $\cl\bigl((A\cap S\cap X)\cup
  D\bigr)=S$, so $$r(A\cap S\cap X)+r(D)\geq r(S)+r(A\cap S\cap X\cap
  D).$$
  Note that $A\cap S\cap X\cap D$ is $A\cap X\cap D$.  Adding
  the three inequalities and simplifying yields the desired
  inequality.

  A similar argument applies to the lattice $C$, for which it suffices
  to consider the antichains $\{A,B,Y\}$, $\{A,W,Y\}$, $\{B,V,Y\}$,
  and $\{V,W,Y\}$.  The last three are listed in MI-orderings.  For
  $\{A,B,Y\}$, apply semimodularity to the pairs $\{A,S\}$, $\{A\cap
  S,B\}$, and $\{A\cap S\cap B,Y\}$; the inclusions $V\subseteq A\cap
  S\subseteq S$ and $T\subseteq A\cap S\cap B\subseteq S$ give
  $\cl\bigl((A\cap S)\cup B\bigr)=E$ and $\cl\bigl((A\cap S\cap B)\cup
  Y\bigr)=S$; add the resulting inequalities and cancel the common
  terms to get inequality (MI) for $\{A,B,Y\}$.
\end{proof}

\begin{figure}[t]
\begin{center}
\epsfxsize 5.75cm \epsffile{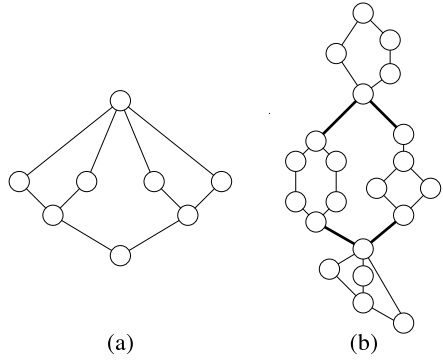} 
\caption{(a) The lattice $L_1*L_2$ where $L_1$ and $L_2$ are Boolean
  lattices on two elements. (b) A lexicographic sum; the indexing
  lattice is a Boolean lattice on two elements.}\label{fig:star}
\end{center}
\end{figure}

We now consider two operations for producing new Tr-lattices.  Given
lattices $L_1$ and $L_2$, let $L_1*L_2$ be the lattice on $(L_1\cup
L_2\cup \{\hat{0},\hat{1}\})-\{\hat{1}_{L_1}, \hat{1}_{L_2}\}$ with
$x\leq y$ if and only if (i) $y=\hat{1}$, or (ii) $x=\hat{0}$, or
(iii) for some $i\in\{1,2\}$, both $x$ and $y$ are in $L_i$ and $x\leq
y$ in $L_i$.  Figure~\ref{fig:star}.a illustrates this operation; note
that the unique four-element antichain in this lattice is not
MI-orderable.

\begin{thm}\label{thm:star}
  If $L_1$ and $L_2$ are Tr-lattices, then so is $L_1*L_2$.
\end{thm}

The proofs of Theorems~\ref{thm:star} and~\ref{thm:lex} are similar,
so we prove only the second result, which concerns lexicographic
sums~\cite[Section~1.10]{tomt}.  Let $L$ be a lattice and let
$\mc{L}=(L_x\,:\,x\in L)$ be a family of lattices that is indexed by
the elements of $L$.  The lexicographic sum $L \oplus \mc{L}$ is
defined on the set $\{(x,a)\,:\, x\in L,\, a\in L_x\}$; the order is
given by $(x,a)\leq (y,b)$ if and only if either (i) $x<y$ in $L$ or
(ii) $x=y$ and $a\leq b$ in $L_x$.  Figure~\ref{fig:star}.b
illustrates this operation.  It is easy to see that $L \oplus \mc{L}$
is not necessarily MI-orderable even if all of the constituent
lattices are.

\begin{thm}\label{thm:lex}
  If $L$ has width at most two and if $\mc{L}=(L_x\,:\,x\in L)$ is a
  family of Tr-lattices, then $L\oplus \mc{L}$ is a Tr-lattice.
\end{thm}

\begin{proof}
  Let $\phi:\mc{Z}(M) \rightarrow L\oplus \mc{L}$ be an isomorphism.
  We must show that any antichain $\mc{A}$ in $\mc{Z}(M)$ satisfies
  inequality (MI).

  For $F\in \mc{Z}(M)$, let $\phi_1(F)$ be the first component of
  $\phi(F)$; thus, $\phi_1(F)\in L$.  For $x\in \phi_1(\mc{A})$, set
  $\mc{A}_x= \{F\,:\, F\in \mc{A},\, \phi_1(F)=x\}$.  Since $L$ has
  width at most two and $\mc{A}$ is an antichain in $\mc{Z}(M)$, there
  are at most two such sets; these sets partition $\mc{A}$.

  For $u\in L$, let $Z_u$ and $E_u$ be the least and greatest flats
  $F\in\mc{Z}(M)$ with $\phi_1(F)=u$.  Thus, if $\phi_1(A)=u$ and
  $\phi_1(B)=v$ with $u\ne v$, then $A\join B=Z_{u\join v}$ and
  $A\meet B=E_{u\meet v}$ by the definition of $L\oplus \mc{L}$.

  Let $x$ be in $\phi_1(\mc{A})$.  Note that $\mc{Z}(M|E_x/Z_x)$ is
  isomorphic to $L_x$, so $M|E_x/Z_x$ is transversal.  Thus, by
  Proposition~\ref{prop:mi},
  $$r_{M|E_x/Z_x}\bigl(\bigcap(\mc{A}_x)-Z_x\bigr) \leq
  \sum_{\mc{F}\subseteq\mc{A}_x} (-1)^{|\mc{F}|+1}
  r_{M|E_x/Z_x}\bigl(\bigcup(\mc{F})-Z_x\bigr),$$
  which gives
  $$r\bigl(\bigcap(\mc{A}_x)\bigr) \leq \sum_{\mc{F}\subseteq\mc{A}_x}
  (-1)^{|\mc{F}|+1} r\bigl(\bigcup(\mc{F})\bigr)$$
  in $M$.  If
  $|\phi_1(\mc{A})|=1$, then the last inequality is the required
  inequality (MI) for $\mc{A}$.  If, instead, $|\phi_1(\mc{A})|=2$,
  let $\phi_1(\mc{A})=\{x,y\}$, so we also have
  $$r\bigl(\bigcap(\mc{A}_y)\bigr) \leq \sum_{\mc{F}\subseteq\mc{A}_y}
  (-1)^{|\mc{F}|+1} r\bigl(\bigcup(\mc{F})\bigr).$$
  The equality
 \begin{align}
   \sum_{\mc{F}\subseteq\mc{A}} (-1)^{|\mc{F}|+1}
   r\bigl(\bigcup(\mc{F})\bigr) = \,& \sum_{\mc{F}\subseteq\mc{A}_x}
   (-1)^{|\mc{F}|+1} r\bigl(\bigcup(\mc{F})\bigr) +
   \sum_{\mc{F}\subseteq\mc{A}_y} (-1)^{|\mc{F}|+1}
   r\bigl(\bigcup(\mc{F})\bigr)    \notag\\
   & \qquad + \sum\limits_{\substack{\mc{F}_x\subseteq\mc{A}_x,\,
       \mc{F}_x\ne\emptyset \\ \mc{F}_y\subseteq\mc{A}_y,\,
       \mc{F}_y\ne\emptyset}} (-1)^{|\mc{F}_x|+|\mc{F}_y|+1}
   r\bigl(\bigcup(\mc{F}_x)\cup \bigcup(\mc{F}_y)\bigr)\notag
 \end{align}
 and that $r(X\cup Y)=r(Z_{x\join y})$ for any $X\in L_x$ and $Y\in
 L_y$ give
  \begin{align}
    \sum_{\mc{F}\subseteq\mc{A}} (-1)^{|\mc{F}|+1}
    r\bigl(\bigcup(\mc{F})\bigr) 
    = \,& \sum_{\mc{F}\subseteq\mc{A}_x} (-1)^{|\mc{F}|+1}
    r\bigl(\bigcup(\mc{F})\bigr) + \sum_{\mc{F}\subseteq\mc{A}_y}
    (-1)^{|\mc{F}|+1}
    r\bigl(\bigcup(\mc{F})\bigr)    \notag\\
    & \qquad -r(Z_{x\join y}) \sum_{\mc{F}_x\subseteq\mc{A}_x,\,
      \mc{F}_x\ne\emptyset} (-1)^{|\mc{F}_x|}
    \sum_{\mc{F}_y\subseteq\mc{A}_y,\, \mc{F}_y\ne\emptyset}
    (-1)^{|\mc{F}_y|}
    \notag\\
    \geq \,& r\bigl(\bigcap(\mc{A}_x)\bigr) +
    r\bigl(\bigcap(\mc{A}_y)\bigr) -r(Z_{x\join y}) \notag\\ \geq \,& 
    r\bigl(\bigcap(\mc{A})\bigr).\notag
  \end{align}
  (The last line uses semimodularity.)  Thus, inequality (MI) holds,
  as needed.
\end{proof}

\section{Open problems}\label{sec:problems}

The results in this paper suggest the following problems.
\begin{enumerate}
\item Is the converse of Theorem~\ref{thm:cover} true? \\ The
  following questions are of interest if the answer is negative.
  \begin{itemize}
  \item[(a)] Find a lattice-theoretic characterization of Tr-lattices,
    perhaps via a recursive description using operations such as those
    in Section~\ref{sec:examples}.
  \item[(b)] Is the converse of Theorem~\ref{thm:cover} true for
    planar lattices?
  \item[(c)] If $\mc{Z}(M)$ has the property of covers in
    Theorem~\ref{thm:cover}, is $M$ a gammoid?
  \item[(d)] Is every sublattice of a Tr-lattice also a Tr-lattice?
    Is this true at least for intervals, or upper intervals?
  \item[(e)] Is the counterpart of Theorem~\ref{thm:ideals} true for
    Tr-lattices?
  \end{itemize}
\item Are there Tr-lattices, or MI-lattices, of all dimensions?
\item Can one capture the minor-closed, dual-closed class of
  transversal matroids described in Theorem~\ref{thm:sublattice} by
  special presentations that such matroids have?  What are the
  excluded minors for this class of matroids?
\item Can one deduce any substantial properties of a matroid $M$ other
  than being a gammoid (or specializations, such as transversal or
  nested) from lattice-theoretic properties of $\mc{Z}(M)$?
\item If $N$ is a minor of $M$ where $\mc{Z}(M)$ is a Tr-lattice, must
  $N$ be transversal?
\item What can we say about $M$ (more particular than transversal)
  when $\mc{Z}(M)$ is a Tr-lattice?
\end{enumerate}

\medskip

\begin{center}
\textsc{Acknowledgements}
\end{center}

I am very grateful for Joseph Kung, whose questions and comments
prompted me to pursue more deeply the implications of
Proposition~\ref{prop:mi}.

\end{document}